\documentclass[12pt,reqno,twoside,a4paper]{amsart}

\usepackage{amssymb, amsmath}

\catcode`@=11 \@mparswitchfalse 
\newcounter{fncount}

\theoremstyle{plain}
\newtheorem{theorem}{Theorem}[section]

\newtheorem{lemma}[theorem]{Lemma}

\newtheorem{fact}[theorem]{Fact}

\newtheorem{pottheorem}[theorem]{False Conjecture}
\newtheorem{restatement}[theorem]{Restatement}

\theoremstyle{definition}
\newtheorem{definition}[theorem]{Definition}

\theoremstyle{remark}
\newtheorem{remark}[theorem]{Remark}

\newcommand{\X}{\mathfrak X}

\newcommand{\Y}{\mathcal Y}
\newcommand{\Z}{\mathcal Z}

\newcommand{\N}{\mathbb N}
\newcommand{\e}{\varepsilon}

\newcommand{\wt}{\widetilde}

\newcommand{\DPbs}{DP-biorthogonal system }
\newcommand{\cobs}{$c_0$-biorthogonal system }
\newcommand{\DPh}{{$\text{DP}_h$ }}
\newcommand{\DPhbs}{{$\text{DP}_h$-biorthogonal system} }
\newcommand{\lnm}{\left\Vert}
\newcommand{\rnm}{\right\Vert}
\newcommand{\lav}{\left\vert}
\newcommand{\rav}{\right\vert}

\newcommand{\lsp}{\left[}
\newcommand{\rsp}{\right]}
\newcommand{\spn}[1]{\text{\rm{sp}}\{ #1 \}}
\renewcommand{\leq}{\leqslant}
\renewcommand{\geq}{\geqslant}

\newcommand{\bio}[2]{\left\{#1,~#2\right\}}
\newcommand{\BIO}[2]{#1~\times~#2}

\newcommand{\ft}{fundamental total}
\newcommand{\nc}{(2+\e)}

\newcommand{\embed}{\hookrightarrow}
\newcommand{\pperp}{\top}  

\newcommand{\deq}{:=}

\newcommand{\hsp}[1]{\hskip #1 pt}

\newcommand{\df}{:=}

\newcommand{\lp}{\left(}
\newcommand{\rp}{\right)}
\newcommand{\linf}{\underline{\textnormal{lim}}}
\newcommand{\lin}[1]{\overset{\underbar{\hspace{#1 pt}}}{\vphantom{,}}} 
\newcommand{\updots}{\overset{\dots}{\vphantom{,}}}

\begin{document}
\baselineskip 18 pt

\title{Biorthogonal Systems in Banach Spaces}

\author{Michael A. Coco}

\address{Department of Mathematics\\
        Lynchburg College\\
         Lynchburg, VA 24501, U.S.A.}
\email{coco@lynchburg.edu}
\thanks{This work is part of the author's PhD dissertation under the direction of
Dr. Maria K. Girardi at the University of South Carolina}

\renewcommand{\subjclassname}{\textup{2000} Mathematics Subject Classification}

\begin{abstract}
We give biorthogonal system characterizations of Banach spaces
that fail the Dunford-Pettis property, contain an isomorphic copy
of $c_0$, or fail the hereditary Dunford-Pettis property.  We
combine this with previous results to show that each infinite
dimensional Banach space has one of three types of biorthogonal
systems.
\end{abstract}

\maketitle


\section{Introduction}
\label{s:intro}
\setcounter{equation}{0}

When we first encounter an arbitrary Banach space, we usually
search for some kind of fundamental structure in the space to make
our understanding of it more complete.  Very often, if a space has
(or fails) a certain property, we can find a fundamental structure
within the space that reflects the property (or failure thereof).
Of course, in this case, we would like to find a strong structure,
like a Schauder basis or finite dimensional decomposition (FDD),
in the space. However, this is not always possible, as even a
separable Banach space need not contain a Schauder basis \cite{E}.
For this reason it is interesting to consider weaker structures
than FDD's and Schauder bases which exist in every separable
Banach space and try to prove that a separable Banach space has a
certain property if and only if there is structure in the space
which reflects the property.

One useful   basis-like structure that has been considered  for a
long time is that of fundamental total biorthogonal system.
Markushevich~\cite{M} showed in 1943 that each separable Banach
space contains a fundamental total biorthogonal system. The main
theorems  of this paper   give a biorthogonal  system
characterization of spaces failing the Dunford-Pettis property and
spaces containing an isomorphic copy of $c_0$.  Combining this
with work already done in the field yields a theorem about the
existence of biorthogonal systems in any given infinite
dimensional Banach space.


\section{Notation and Motivation}
\label{s:notation}

Throughout this paper, $\X$
denotes  an arbitrary
(infinite-dimensional real) Banach space.
If $\X$ is a Banach space,
then
$\X^*$  is its topological dual space,
$B(\X)$ is its  (closed) unit ball, and
 $S(\X)$ is its unit sphere.
If $X$ is a subset of $\X$, then
$\spn  X $ is the linear span of $X$ while
$\lsp X\rsp$  is
the closed linear span of~$X$.
The Kronecker delta  $\delta_{nm}$ takes the value 1 when
$n=m$ and 0 when $n\neq m$.

\begin{definition} For a subset $X$ of $\X$
and a subset $Z$  of $\X^*$:
\begin{enumerate}
\item the \textit{annihilator} of $X$ is $X^\perp=\{x^*\in\X^*:x^*(x)=0 \text{ for all } x\in X\}$,
\item the \textit{preannihilator} of $Z$ is $Z^\pperp=\{x\in\X: x^*(x)=0 \text{ for all } x^*\in Z\}$,
\item
$X$  is  \textit{fundamental} if $\lsp X \rsp = \X$,
or,  equivalently, $X^\perp=\{ 0 \}$,
\item
$Z$  is  \textit{total} if the  weak$^*$-closure
of $\spn{Z}$ is $\X^*$,
or, equivalently,  \newline$Z^\pperp=\{ 0 \}$,
\item
for  a fixed $\tau \geq 1$,
$Z$ \textit{$\tau$-norms} $X$
(or $X$ is $\tau$-normed by $Z$)  if
\[
\lnm x \rnm~ \leq~ \tau \, \sup_{z\in Z\setminus \{0\}}
\frac{z(x)}{\lnm z \rnm}
\]
for each ~$x\in X$,
\item
$Z$   \textit{norms} $X$ if $Z$   $1$-norms  $X$.
\end{enumerate}
\end{definition}
\noindent It is easy to see that if $Z$ $\tau$-norms $\X$ for a $\tau \geq 1$
then $Z$ is total.

\begin{definition} A system $\bio{x_n}{x_n^*}_{n=1}^{\infty}$
in $\BIO{X}{Z}$ is
\begin{enumerate}
\item
a   \textit{biorthogonal system}
if $x_n^*(x_m) = \delta_{nm}$,
\item
\textit{$M$-bounded}
if  \,
$\{ x_n \}$ and $\{ x^*_n \}$ are bounded
and \newline $ \sup_n \lnm x_n \rnm \, \lnm x_n^* \rnm
\leq M$,
\item
\textit{bounded} if it is $M$-bounded for some (finite) $M$,
\item
\textit{fundamental} if $\{ x_n  \}$ is fundamental,
\item
\textit{total} if   $\{x_n^*\}$ is  total.
\end{enumerate}
\end{definition}

A sequence $\{x_n\}_{n=1}^\infty$ in a Banach space $\X$ is called
\textit{semi-normalized} if there are constants $0<\alpha\leq\beta<\infty$
such that $\alpha\leq\|x_n\|\leq\beta$ for each $n\in\N$. Recall that
$\{x_n\}_{n=1}^\infty$ is \textit{a basic sequence} if each $x_n$ is non-zero
and there exists a finite constant $K>0$ such that
\begin{equation}\label{bas}\left\|\sum_{j=1}^m a_jx_j\right\|\leq K\left\|\sum_{j=1}^n a_jx_j\right\|\end{equation}
for all choices $\{a_j\}_{j\in\N}$ and any integers $m<n$.  When
this is the case, the smallest $K$ for which \eqref{bas} holds is
called the \textit{basis constant} of $\{x_n\}_{n=1}^\infty$ and
there exists a biorthogonal system $\bio{x_n}{x_n^*}$ in
$\BIO{\X}{\X^*}$ such that $\|x_n^*\|\leq\frac{2K}{\|x_n\|}$.

\noindent Operators between Banach spaces are assumed to be bounded and linear.
All notation and terminology, not otherwise explained,
are as in~\cite{DU} or~\cite{LT1}.

Our motivation begins with the following structure theorem of E. Odell \cite{O}:

\begin{theorem}
\label{t:Bs}
Every infinite dimensional Banach space contains a subspace isomorphic to $c_0$, a
subspace isomorphic to $\ell_1$ or a subspace that fails the Dunford-Pettis property.
\end{theorem}

\noindent Our goal is to find a biorthogonal system version of this theorem in
which the conditions imposed on the biorthogonal systems directly reflect the
property they characterize.  Luckily, some of the work, the $\ell_1$ case, has
already been done for us.  In fact, our results are inspired by this previous work.
In 2000, S.J.~ Dilworth, M.~ Girardi, and W.B.~Johnson characterized spaces containing
isomorphic copies of $\ell_1$ using biorthogonal systems.

\noindent
\begin{theorem}\cite{DGJ}
\label{t:elle}
The following statements  are equivalent.
\begin{enumerate}
\item $\ell_1 \embed \X$.
\item There is a bounded $wc_0^*$-stable
             biorthogonal system
             in $\BIO{\X}{\X^*}$.
\end{enumerate}
And in the case that $\X$ is separable:
\begin{enumerate}
\item[(3)] There is a bounded \ft\  $wc_0^*$-stable
             biorthogonal system
             $\bio{x_n}{x_n^*}$ in $\BIO{\X}{\X^*}$.
\end{enumerate}
Furthermore for each $\e>0$:
if {\rm (2)} holds
then the system can be taken to be $(1+\e)$-bounded;
if {\rm (3)}
holds then the system can be taken to be $[(1+\sqrt2)  + \e]$-bounded
and so that  $\lsp x_n^* \rsp $  $\nc$-norms   $\X$.
\end{theorem}

Recall that $\bio{x_n}{x_n^*}$ is a \textit{$wc_0^*$-stable biorthogonal system}
if, for each isomorphic
embedding $T$ of $\X $ into some  $ \Y$,
there exists a lifting $\{ y_n^* \}$ of $\{ x_n^* \}$
(i.e., $T^* y_n^* = x_n^*$ for each $n$)
such that
$\{ y_n^* \}$ is a semi-normalized weakly-null sequence in $\Y^*$
(or equivalently, such that   $\bio{Tx_n}{y_n^*}$ in $\BIO{\Y}{\Y^*}$
is a $wc_0^*$-biorthogonal system).

\noindent They also characterized Banach spaces that have
\textit{Schur property} (i.e. weak  and  strong sequential
convergence in $\X$ coincide) via Biorthogonal systems.  In the
next section we will discuss the Dunford-Pettis property.  Recall
that the Schur property is related to the Dunford-Pettis property
and embeddings of $\ell_1$ in the following way:
(cf.~\cite[p.~23]{D2}) $\X^*$ fails the Schur property if and only
if $\X$ fails the Dunford-Pettis property or $\ell_1 \embed \X$.
This fact provides a link between the above results and the
results of the next section that characterize failure of the
Dunford-Pettis Property.


\section{Spaces Failing The Dunford-Pettis Property}
\label{s:dp}

Recall that a Banach space $\X$ has the Dunford-Pettis property (DP) if
whenever $\{x_n\}_n\subset\X$ and $\{x_n^*\}_n\subset\X^*$ are weakly null
sequences, we have $\lim_{n\to\infty}x_n^*(x_n)=0$.  We refer the reader to
the excellent survey article \cite{D2} for a complete treatment of all things Dunford-Pettis.  Further results and additional open questions can be found in \cite{CG}.

Now suppose $\X$ is a Banach space that fails the  Dunford-Pettis property.
Then  there exists
a weakly null sequence $\left\{ w_k \right\}_{k\in\N} $  in $\X$
and  a weakly null sequence $\left\{ w_k^* \right\}_{k\in\N} $  in $\X^*$
  such that $\lim_{k\to\infty} \lav w_k^* \lp w_k \rp\rav  \neq 0$.
We may assume, without loss of generality, that there exists $\delta>0$
such that $ w_k^*(w_k) >\delta$ for each $k\in\N$.  If this is not the case we can
pass to a suitable subsequence and adjust signs.
Now $\left\{ w_k \right\}_{k\in\N} $ and
$\left\{ w_k^* \right\}_{k\in\N} $  are semi-normalized so we may renormalize
if necessary to get that for each $k\in\N$:
\begin{enumerate}
\item $w_k \in   S(\X)$,
\item $w^*_k\lp w_k \rp = 1$,
\item $1\leq \lnm w_k^*\rnm \leq M$ for some constant $M$.
\end{enumerate}
This leads to the following definition.
\begin{definition}
\label{d:mdpp}
Let $M\geq 1$.
$\X$ fails the \textit{$M$-Dunford-Pettis property}
provided
there is a weakly null sequence $\{w_k\}_k $ from $   S(\X)$
and a weakly null sequence $\{w_k^*\}_k $ from $ \X^*$ such that
$w^*_k\lp w_k \rp = 1$ and $1\leq \lnm w_k^*\rnm \leq M$ for each $k\in\N$.
\end{definition}
Note that clearly $\X$ fails $M$-DP for some $M$ if and only if $\X$ fails DP.
We only bother to define it here to make the statement of Theorem \ref{t:dpit} a bit clearer.

\begin{definition} A biorthogonal system $\bio{x_n}{x_n^*}$ in $\BIO{\X}{\X^*}$
is called a  \textit{\DPbs}
if $\{x_n  \}$ and $\{x_n^* \}$  are semi-normalized weakly-null  sequences.
\end{definition}

\begin{theorem}
\label{t:dpit}
The following statements  are equivalent.
\begin{enumerate}
\item $\X$  fails the Dunford-Pettis property.
\item There is a  bounded  \DPbs
             in $\BIO{\X}{\X^*}$.
\end{enumerate}
And in the case that $\X$ is separable:
\begin{enumerate}
\item[(3)] There is a  bounded \ft\ \DPbs  \\
 $\bio{x_n}{x_n^*}$ in $\BIO{\X}{\X^*}$.
\end{enumerate}
Furthermore, for an $\X$  failing the $M$-Dunford-Pettis property,
  for each $\e>0$:
if {\rm (2)} holds then the system can be taken to be $(M+\e)$-bounded;
if {\rm (3)} holds then the system
can be taken to be $[M(1+\sqrt2)^2 + \e]$-bounded
and so that  $\lsp x_n^* \rsp $ norms $\X$.
\end{theorem}
It is clear that (2)  implies (1) as well as (3) implies (1).
That (1) implies (2) follows from
Theorem~\ref{t:dp}.
That (1) implies (3) in the separable case follows from
Theorem~\ref{t:ftdp}.

The following well-known basic fact will be used.
\begin{fact}
\label{f:bf}
Let $\X_0$ be a finite codimensional subspace of $\X$
and $\left\{ x_n \right\}_{n\in\N}$ be a weakly null sequence in ~$\X$.
Then
\begin{equation*}
d(x_n, \X_0) \df \inf_{x_0\in\X_0} \lnm x_n - x_0 \rnm
~\xrightarrow{~n\to\infty~}~ 0~.
\end{equation*}
Thus, if $\{x_n\}_n$ is semi-normalized and $\e >0$,
 there exists $n_\e$ and    $\widetilde{x}_{n_\e}\in\X_0$
with $\lnm   x_{n_\e} - \widetilde{x}_{n_\e}  \rnm < \e$
and $\lnm  {x}_{n_\e} \rnm = \lnm \widetilde{x} _{n_\e} \rnm$.
\end{fact}

We can now give a quantitative proof that (1) implies (2) in Theorem \ref{t:dpit}.

\begin{theorem}
\label{t:dp}
Let $\X$ fail  the $M$-Dunford-Pettis property
 and $\e>0$.  Then there is a biorthogonal system
$\bio{x_n}{x_n^*}_{n=1}^\infty$
in $\BIO{\X}{\X^*}$  such that:
\begin{enumerate}
\item $\left\{ x_n \right\}_{n=1}^\infty$
and   $\left\{ x_n^* \right\}_{n=1}^\infty$
are weakly null
\item $\lnm x_n \rnm = 1$ for each $n\in\N$
\item  $1\leq    \lnm x_n^* \rnm  \leq M +\e$ for each $n\in\N$
\item $\left\{ x_n \right\}_{n=1}^\infty$ is a basic sequence.
\end{enumerate}
\end{theorem}

\begin{proof}
Since $\X$ fails  the $M$-Dunford-Pettis property  there exist
sequences $\{ w_k \}_{k\in\N} $ and    $\left\{ w_k^*
\right\}_{k\in\N} $   as in Definition~\ref{d:mdpp}. Without loss
of generality  (pass to a subsequence)   $\{ w_k \}_{k\in\N} $ is
a basic sequence.

Let $\left\{ \e_n \right\}_{n\in\N}$ be a decreasing sequence of
positive numbers with $\e_1 < \frac{\e}{2(M+\e)}$ and
$\sum_{n\in\N} \e_n < \frac{1}{2K}$ where $K$ is the basis
constant of   $\left\{ w_k \right\}_{k\in\N} $. We will construct
a   system $\bio{x_n}{x_n^*}_{n=1}^\infty$ in $\BIO{\X}{\X^*}$ and
an increasing sequence $\{ k_n \}_{n\geq 1}$ of integers  such
that
\begin{enumerate}
\item[(a)] $\bio{x_n}{x_n^*}_{n=1}^\infty$ is biorthogonal
\item[(b)] $\lnm x_n \rnm = 1$
\hsp{12} for each $n\in\N$
\item[(c)] $ 1 \leq \lnm x_n^*  \rnm~ \leq ~\frac{M}{1-2\e_n}$
\hsp{12} for each $n\in\N$
\item[(d)] $ \lnm x_n - w_{k_n} \rnm~ \leq ~\frac{\e_n}{M}$
\hsp{12}  for each $n\in\N$
\item[(e)] $ \lnm x_n^* - w_{k_n}^* \rnm ~\leq ~
  \e_n +   \frac{2M\e_n}{1-2\e_n} $
\hsp{12}   for each $n\in\N$.
\end{enumerate}
\noindent Conditions (d) and (e) will give us (1): for $x^*\in\X^*$
$$|x^*(x_n)|\leq\|x^*\|\|x_n-w_{k_n}\|+\left|x^*(w_{k_n})\right|\to 0$$
so $\{x_n\}_n$ is weakly null and simlilarly for $\{x_n^*\}_n$.

\noindent Condition (c) gives us (3):
$$1\leq\|x_n^*\|\leq\frac{M}{1-2\e_n}\leq\frac{M}{1-2(\frac{\e}{2(M+\e)})}=M+\e.$$

\noindent Condition (d) gives us (4): we have
$$\sum_n\|w_{k_n}-x_n\|\leq\sum_n\e_n<\frac{1}{2K}.$$
Then $\{x_n\}_n$ is basic (and equivalent to $\{w_{n_k}\}_k$).

Now we construct $\bio{x_n}{x_n^*}_{n=1}^\infty$ by induction. To
start,  let $k_1 = 1$ and $x_1  = w_1 $ and     $x_1^* = w_1^*$.
Fix $n > 1$ and assume that a system  $\bio{x_j}{x_j^*}_{j < n}$,
along with a sequence $\{ k_j \}_{j < n}$, have been constructed
to satisfy the above conditions. Let
$$\X_n =\lsp x_j^*\rsp_{j < n}^\pperp
\text{\quad and \quad}
 \Z_n =\lsp x_j\rsp_{j < n}^\perp \ .
$$
Using Fact~\ref{f:bf},
find $k_n > k_{n-1}$ and $x_n\in\X_n$ and $z_n^*\in\Z_n$  so that
\begin{equation*}
d\lp w_{k_n}, \X_n \rp
~\leq~
\lnm w_{k_n} - x_n \rnm
~<~ \frac{\e_n}{M}
\hsp{10}\textrm{and}\hsp{10}
 d\lp w_{k_n}^*, \Z_n \rp
~\leq~
\lnm w_{k_n}^* - z_n^* \rnm
~<~    \e_n
 \end{equation*}
 with
\begin{equation*}
\lnm x_n \rnm = 1
\hsp{15}\textrm{and}\hsp{15}
1 \leq \lnm z^*_n \rnm \leq M ~.
\end{equation*}
Note that
\begin{align*}
\lav z^*_n \lp x_n \rp - w^*_{k_n} \lp  w_{k_n}  \rp \rav &~=~
\lav z_n^* \lp x_n - w_{k_n} \rp
- \lp w^*_{k_n} - z_n^* \rp \lp w_{k_n} \rp \rav \\
&<~ M \ \frac{\e_n}{M} ~+~ \e_n ~=~ 2\ \e_n ~,
\end{align*}
and so $1-2\e_n < z^*_n \lp x_n \rp  < 1+2 \e_n.$  Let
\begin{equation*}
x_n^* ~\df~ \frac{z_n^*}{z^*_n\lp x_n \rp} ~.
\end{equation*}
 Thus conditions (a) and (c) hold.  As for condition (e):
\begin{align*}
\|x_n^*-w_{k_n}^*\|&\leq\|w_{k_n}^*-z_n^*\|+\|z_n^*-\frac{z_n^*}{z_n^*(x_n)}\|\\
&\leq \e_n+\frac{1}{z_n^*(x_n)}|z_n^*(x_n)-1|\|z_n^*\|\\
&\leq\e_n +\frac{2\e_n}{1-2\e}M.
\end{align*}
\end{proof}

The construction of  fundamental total
biorthogonal systems  in the
proofs of~(1) implies~(3) in  Theorem~\ref{t:dpit}
 and Theorem \ref{t:c0it}
use  the Haar matrices,
which are summarized below.

\begin{remark}
\label{r:hm}
Fix $m \geq 0$ and consider the
$2^m$-dimensional Hilbert space $\ell_2^{2^m}$,
along with its unit vector basis $\{ e^2_j \}_{j=1}^{2^m}$.

The   Haar basis $\{ h^m_j \}_{j=1}^{2^m}$
of $\ell_2^{2^m}$ can be described as follows.  For $0 \leq n \leq m$
and $1 \leq k \leq 2^n$  let
\[
 I^n_k ~=~ \left\{ j \in \mathbb N \hsp{3}\colon \hsp{3}
2^{m-n} \, (k-1) \hsp{3}<\hsp{3}
j \hsp{3}\leq\hsp{3}
2^{m-n} \,  k  \right\}\ .
\]
Thus
\begin{gather*}
I^0_1 = \left\{ 1, 2 \hsp{2},\hsp{2}  \ldots \hsp{2},\hsp{2} 2^m \right\} \\
I^1_1  = \left\{ 1, 2 \hsp{2},\hsp{2}\ldots \hsp{2},\hsp{2}
          2^{m-1} \right\}
\textrm{\hsp{10}and\hsp{10}}
I^1_1  = \left\{ 1 + 2^{m-1}  \hsp{2},\hsp{2} \ldots \hsp{2},\hsp{2}
            2^m \right\} \ .
\end{gather*}
 In general,
the collection $\{ I^n_k \}_{k=1}^{2^n}$
of sets along the  $n^{\text{th\/}}$-level (disjointly) partitions
the set
$\{ 1, 2, \ldots , 2^m \}$ into $2^n$ sets, each containing
$2^{m-n}$ consecutive integers,
and $I^n_k$ is the disjoint union $I^n_k = I^{n+1}_{2k-1}  \cup I^{n+1}_{2k} $.
Now let
\[
h^m_1 ~=~ 2^{\frac{-m}{2}} ~ \sum_{j\in I^0_1} e^2_j
\]
and, for $0 \leq n < m$ and $1 \leq k \leq 2^n$,
let $h^m_{2^n + k}$ be supported on $I^n_k$ as
\[
h^m_{2^n + k} ~=~ 2^{\frac{n-m}{2}} ~
\left[ \sum_{j\in I^{n+1}_{2k-1}} e^2_j
~-~ \sum_{j\in I^{n+1}_{2k}} e^2_j  \right] \ .
\]
Note that   $\{ h^m_j \}_{j=1}^{2^m}$
forms an orthonormal basis for~$\ell_2^{2^m}$.

Let $H_m = \left( a^m_{ij} \right)$ be the $2^m \times 2^m$ Haar
matrix  that transforms the unit vector basis  of $\ell_2^{2^m}$
onto the Haar basis; thus, the $j^{\text{th\,}}$ column vector of
$H_m$ is just  $h^m_{j}$ and so $H_m$ is a unitary matrix. For
example, for $m=2$ we have
$$
H_2 ~=~  \bmatrix
2^{-1}  &  +2^{-1}  & +2^{-1/2}  & 0 \\
2^{-1}  &  +2^{-1}  & -2^{-1/2}  & 0 \\
2^{-1}  &  -2^{-1} &  0          & +2^{-1/2}  \\
2^{-1}  &  -2^{-1} &  0          & -2^{-1/2}
\endbmatrix ~ .
$$

Now if $\bio{z_j}{z_j^*}_{j=1}^{2^m}$  is a biorthogonal sequence
in $\BIO{\X}{\X^*}$ and $\bio{x_i}{x_i^*}_{i=1}^{2^m}$ is such
that
\begin{align*}
H_m ~ \bmatrix  z_1 \\ \vdots \\  z_{2^m} \endbmatrix
~=~ \bmatrix  x_1 \\ \vdots \\  x_{2^m} \endbmatrix
\hskip 40 pt &\text{and}\hskip 40 pt
H_m ~ \bmatrix  z_1^* \\ \vdots \\  z_{2^m}^* \endbmatrix
~=~ \bmatrix   x_1^* \\ \vdots \\   x_{2^m}^* \endbmatrix \ , \\
\intertext{then} x_i \deq \sum_{j=1}^{2^m} a^m_{ij}  z_{j} \hskip
40 pt &\text{and}\hskip 40 pt x_i^*  \deq \sum_{j=1}^{2^m}
a^m_{ij} z_{j}^*   \ .
\end{align*}
It is not hard to see that since $H_m $ is a unitary  matrix,
\begin{enumerate}
\item[(1)] $x_i^* (x_j) = \delta_{ij} $ \item[(2)] $\lsp x_i
\rsp_{i=1}^{2^m}
        = \lsp   z_j \rsp_{j=1}^{2^m} $
\item[(3)]  $\lsp x_i^* \rsp_{i=1}^{2^m}
       = \lsp   z_j^* \rsp_{j=1}^{2^m} $.
\end{enumerate}
Note that, for each $1\leq i\leq 2^m$,
\begin{enumerate}
\item[(4)]  $ a_{i1}^m = 2^{{-m}/{2}} $
\end{enumerate}
and
\begin{enumerate}
\item[(5)]  $\sum _{j=2}^{2^m} \lav a^m_{ij}\rav
            ~=~ \lp 1 + \sqrt 2 \rp \
                \lp 1 - 2^{\frac{-m}{2}}  \rp
              ~  \ ^{\ _{m \to \infty}} \hskip - 8 pt \nearrow
              ~ 1 + \sqrt 2 $.
\end{enumerate}
It follows that
\begin{enumerate}
\item[(6)] $\lnm   x_i \rnm \ \ ~ \leq ~
     2^{{-m}/{2}} \lnm  z_{1}  \rnm \   ~+~
      \left( 1 + \sqrt 2  \right)
           \max_{ 1 < j \leq 2^m }   \lnm   z_j  \rnm$
\item[(7)]  $\lnm  x_i^* \rnm ~ \leq ~
       2^{{-m}/{2}} \lnm  z_{1}^* \rnm   ~+~
      \left( 1 + \sqrt 2  \right)
           \max_{ 1 < j \leq 2^m }   \lnm   z_j^* \rnm$
\item[(8)] for each $x^{*} \in \X^{*}$  \newline
  $\lav x^{*} (x_i) \rav    ~\leq~  $
    $  2^{{-m}/{2}}  \lav x^{*} \lp z_{1}  \rp \rav  ~+~
           \left( 1 + \sqrt 2  \right)
          \max_{ 1 < j \leq  2^m }
         \lav  x^{*} \left(   z_j  \right) \rav $
\item[(9)] for each $x^{**} \in \X^{**}$  \newline
  $\lav x^{**} (x_i^*) \rav    ~\leq~  $
    $  2^{{-m}/{2}}  \lav x^{**} \lp z_{1}^*  \rp \rav  ~+~
           \left( 1 + \sqrt 2  \right)
          \max_{ 1 < j \leq  2^m }
         \lav  x^{**} \left(   z_j^*  \right) \rav $ .
\end{enumerate}
\end{remark}
\noindent
The following notation will (hopefully) simplify the
proofs of~Theorem~\ref{t:ftdp} and Theorem \ref{t:ftc0}.
\begin{definition}
\label{d:blk} A sequence $\{ J_k \}_{k=1}^\infty$ of
subsets of $\mathbb N$ is a  \textit{blocking}
of $\mathbb N$ if $\mathbb N$ is the
disjoint union $\cup_{k=1}^\infty J_k$ and
$$
   \max J_k < \min J_{k+1}
$$
for each $k\in \mathbb N$.  Given a blocking  $\{ J_k \}_{k=1}^\infty$
of $\mathbb N$, let $J_0 = \emptyset $ and
\begin{alignat*}{2}
J_k^p   ~&\deq~  \bigcup\limits_{0\leq j<k} J_j
&\hsp{30}\textrm{,}\hsp{30}
J_k^o  ~&\deq~ J_k\setminus\{ \text{the first element in } J_k\}  \\
J_k^{po}  ~&\deq~  \bigcup\limits_{0\leq j<k}  J_j^o
&\hsp{30}\textrm{,}\hsp{30}
\mathbb N^o ~&\deq~ \bigcup\limits_{k=1}^\infty  J_k^o
\end{alignat*}
for each $k \in \mathbb N$.  Pictorially one has:
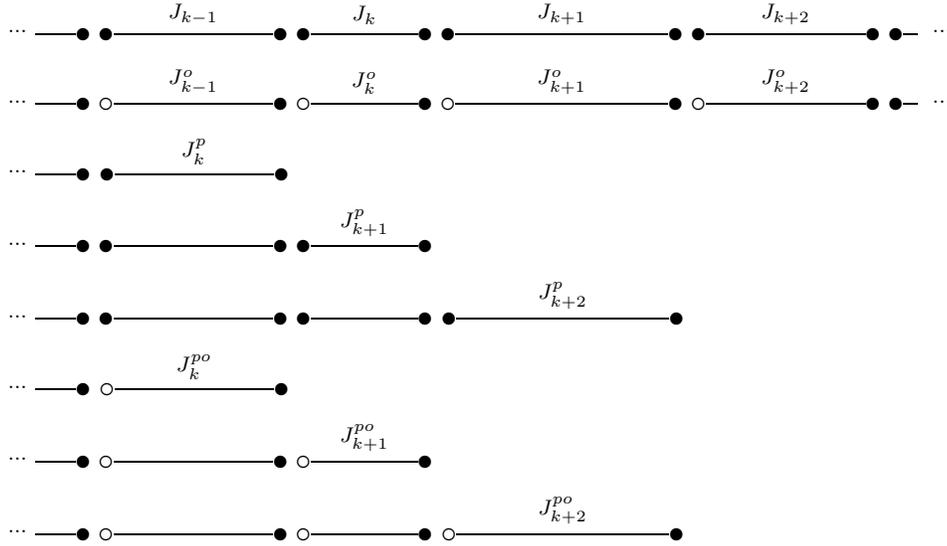
\begin{align*}     
\updots\; \lin{15}\bullet &\overset{J_{k-1}}{\bullet\lin{60}\bullet}\overset{J_{k}}{\bullet\lin{40}\bullet} \overset{J_{k+1}}{\bullet\lin{80}\bullet} \overset{J_{k+2}}{\bullet\lin{60}\bullet}\overset{}{\bullet\lin{5}}\; \updots\\
\updots\; \lin{15}\bullet &\overset{J_{k-1}^o}{\circ\lin{60}\bullet}\overset{J_{k}^o}{\circ\lin{40}\bullet} \overset{J_{k+1}^o}{\circ\lin{80}\bullet} \overset{J_{k+2}^o}{\circ\lin{60}\bullet}\overset{}{\bullet\lin{5}}\; \updots\\
\updots\; \lin{15}\bullet &\hskip3pt\overset{J_{k}^p}{\bullet\lin{60}\bullet}\\
\updots\; \lin{15}\bullet &\overset{}{\bullet\lin{60}\bullet}\overset{J_{k+1}^p}{\bullet\lin{40}\bullet}\\
\updots\; \lin{15}\bullet &\overset{}{\bullet\lin{60}\bullet}\overset{}{\bullet\lin{40}\bullet}\hskip3pt\overset{J_{k+2}^p}{\bullet\lin{80}\bullet}\\
\updots\; \lin{15}\bullet &\hskip3pt\overset{J_{k}^{po}}{\circ\lin{60}\bullet}\\
\updots\; \lin{15}\bullet &\overset{}{\circ\lin{60}\bullet}\overset{J_{k+1}^{po}}{\circ\lin{40}\bullet}\\
\updots\; \lin{15}\bullet &\overset{}{\circ\lin{60}\bullet}\overset{}{\circ\lin{40}\bullet}\hskip3pt\overset{J_{k+2}^{po}}{\circ\lin{80}\bullet}\\
\end{align*}
\end{definition}

It follows from the next theorem that (1) implies (3) for
separable $\X$ in Theorem \ref{t:dpit}.

\begin{theorem}
\label{t:ftdp}
Let $\X$ fail the $M$-Dunford-Pettis property and $\e > 0$.  If
$\bio{ a_n }{ b_n^*}_{n\in\mathbb N} \subset
\BIO{\X}{\X^*}$.
then there exists a
 $[M(1+\sqrt{2})^2+\e]$-bounded
\DPbs $\bio{x_n}{x_n^*}$ in $\BIO{\X}{\X^*}$ such that
$\lsp a_n \rsp_{n\in\N} \subset \lsp x_n \rsp_{n\in\N}  $
and $\lsp b_n^*  \rsp_{n\in\N} \subset \lsp x_n^* \rsp_{n\in\N}  $.
\end{theorem}

\begin{proof}
Without loss of generality,
$[ a_n ]_{n\in\mathbb N}$ and  $[ b^*_n ]_{n\in\mathbb N}$
are each infinite dimensional.  Since $\X$ fails the
$M$-Dunford-Pettis property,
by Theorem~\ref{t:dp},
there is a biorthogonal system $\bio{w_n}{w_n^*}$
in $\BIO{\X}{\X^*}$
with both $\{w_n\}_n$ and $\{w_n^*\}_n$ weakly null,
$\|w_n\|=1$, and $1\leq\|w_n^*\|\leq M+\e$.
Fix a sequence $\{ \delta_k \}_{k=1}^\infty$
of positive numbers decreasing to zero with $\delta_1<\frac{1}{2}$ and
\begin{equation}
\label{e:Mde}
\frac{M+\e}{\lp 1+2\e\rp M}
~<~
1-2\delta_1 ~.
\end{equation}

It suffices to find a system
$\bio{x_n}{x_n^*}_{n=1}^\infty$ in $\BIO{\X}{\X^*}$
along with
(following the terminology in Definition~\ref{d:blk})
a blocking $\{ J_k \}_{k=1}^\infty$ of~$\mathbb N$
and
an increasing sequence $\{ i_n\}_{n\in \mathbb N^o}$ from $\mathbb N$,
satisfying
\begin{enumerate}
\item[\rm{(1)}]
 $x_m^*(x_n) = \delta_{mn}$
\item[\rm{(2)}]
 $\lnm x_n \rnm ~\leq~ \left(1+\sqrt{2} \right) + \e$
\item[\rm{(3)}]
 $\lnm x_n^* \rnm ~\leq~  \lp 1+2\e\rp M\left(1+\sqrt{2} \right) + \e $
\item[\rm{(4)}]
for each $x^*\in S\left(\X^*\right)$, if $n\in J_k$, then
\newline
$\lav x^{*}\left(x_n\right) \rav ~\leq~
\delta_k + \left( 1+\sqrt2\right) \, \max_{j \in J^o_k}
\left(~ \lav x^{*}\left( w_{i_j} \right) \rav
      +  \delta_k~\right)$
\item[\rm{(5)}]
 for each $x^{**}\in S\left(X^{**}\right)$, if $n\in J_k$ then
\newline
$\lav x^{**}\left(x_n^*\right) \rav ~\leq~\delta_k\lp\frac{4+2M}{1-2\delta_k}\rp+\lp 1+\sqrt{2}\rp \max_{j\in J_k^o}\lav x^{**}\lp w_{i_j}^*\rp\rav$

\item[\rm{(6)}]
 $\lsp a_n \rsp_{n=1}^{\infty}
       \subset \lsp x_n \rsp_{n=1}^{\infty}$
\item[\rm{(7)}]
  $\lsp b^*_n \rsp_{n=1}^{\infty}
       \subset \lsp x_n^* \rsp_{n=1}^{\infty} $   \ .
\end{enumerate}

The construction will inductively produce
blocks $\bio{x_n}{x_n^*}_{n\in J_k}$.
Let $x_0$ and $x^*_0$ be the zero vectors.
Fix $k \geq 1$.
Assume that $\{ J_j \}_{0\leq j<k}$
along with
$\bio{x_n}{x_n^*}_{n \in J_k^p}$
and  $\{ i_n \}_{n \in J_k^{po}}$ have been constructed to satisfy conditions~(1) through (5).
Now to construct $J_k$
along with $\bio{x_n}{x_n^*}_{n \in J_k }$
and   $\{ i_n \}_{n  \in J_k^o }$.

Let
$$
 \mathcal P_k ~\deq~
    \left[ x^*_n \right]_{n\in J_k^p }^\pperp
    \text{\hskip 20 pt and \hskip 20 pt}
 \mathcal  Q_k ~\deq~
    \left[ x_n \right]_{n \in J_k^p}^\perp
$$
and
$$
    n_k = \max  J^p_k \ .
$$
The idea is to    find a biorthogonal system
$\bio{z_n}{z_n^*}_{n\in J_k}$ in $\BIO{\mathcal  P_k}{\mathcal  Q_k}$
by first finding  just one pair $\{z_{1+n_k}, z^*_{1+n_k} \}$ which helps
guarantee condition~(6) if $k$ is odd  and
condition~(7) if $k$   even;
however,   $\{z_{1+n_k}, z^*_{1+n_k} \}$
would not necessarily satisfy conditions~(2) through~(5)
and so $J^o_k$ and
$$
    \bio{z_n}{z_n^*}_{n \in J_k^o} \ ,
$$
and $\{ i_n \}_{n\in J^o_k}$ are constructed and then the appropriate Haar
matrix is applied to
$\bio{z_n}{z_n^*}_{n\in J_k}$ to produce $\bio{x_n}{x_n^*}_{n\in J_k}$
so that
$$
  \bio{x_n}{x_n^*}_{n\in  J^p_k \, \cup J_k }
$$
with $\{ i_n \}_{n\in J^{po}_k \cup J^o_k}$
satisfy conditions~(1) through (5).

$\bio{z_{1+n_k}}{z^*_{1+n_k}}$ is constructed
by a standard Gram-Schmidt biorthogonal procedure.
If $k$ is odd, start  in~$\X$. Let
$$
h_k ~=~ \min \left\{ h \colon
     a_h \not\in \lsp x_n \rsp_{n \leq n_k} \right\} \ .
$$
Set
\begin{align*}
z_{1+n_k}
~&=~
a_{h_k} ~-~
   \sum_{n \leq n_k} x_n^*(a_{h_k}) x_n \ ,  \\
\intertext{and for any $y^*_{1+n_k}$ in $\X^*$ such
that   $y^*_{1+n_k}(z_{1+n_k}) \neq 0$,}
z^*_{1+n_k}
~&=~
\frac{
y^*_{1+n_k}  ~-~
   \sum_{n \leq n_k}  y^*_{1+n_k}  (x_n) x_n^*}
{  y^*_{1+n_k}  (z_{1+n_k})} \ .
\end{align*}
If $k$ is even, start  in~$\X^*$. Let
$$
h_{k} ~=~ \min \left\{ h \colon
     b^*_h \not\in \lsp x_n^* \rsp_{n\leq n_k} \right\} \ .
$$
Set
\begin{align*}
z_{1+n_k}^*
~&=~
b^*_{h_k} ~-~
   \sum_{n \leq n_k} b^*_{h_k} (x_n) x_n^* \ ,  \\
\intertext{and, for any $y_{1+n_k}$ in $\X$  such
that   $ z_{1+n_k}^* (y_{1+n_k})  \neq 0$,}
z_{1+n_k}
~&=~
\frac{
y_{1+n_k}  ~-~
   \sum_{n\leq n_k}  x^*_n ( y_{1+n_k} )  x_n}
{ z_{1+n_k}^* (y_{1+n_k}) } \ .
\end{align*}
Clearly  $z_{1+n_k}^* \left( z_{1+n_k} \right) = 1$
and
$$
  z_{1+n_k} ~ \in ~ \mathcal  P_k
  \text{\hskip 25 pt and \hskip 25 pt}
    z_{1+n_k}^* ~ \in ~ \mathcal  Q_k   \ .
$$

Find a natural number  $m_k$ larger than one  so that
$$
2^{{-m_k}/{2}} \
\max\left(~ \lnm z_{1+n_k} \rnm ~,~ \lnm z_{1+n_k}^* \rnm ~\right)
~<~
\min\left(~ \e ~,~ \delta_k ~\right)
$$
and let
$$
  J_k \deq \{ 1 + n_k , \ldots ,  2^{m_k} + n_k \}
\text{\hskip 10 pt and  so \hskip 10 pt}
  J_k^o \deq \{ 2 + n_k , \ldots , 2^{m_k} + n_k \} \ .
$$

Let
$$
\wt{\mathcal  P}_k ~\deq~  \mathcal  P_k \cap \left[ z^*_{1+n_k}\right]^\pperp
\text{\hskip 25 pt and \hskip 25 pt}
\wt{\mathcal  Q}_k ~\deq~ \mathcal  Q_k \cap \left[ z_{1+n_k}\right]^\perp \ .
$$
The next step is to   find a biorthogonal system
$\bio{z_n}{z_n^*}_{n\in J_k^o}$
along with $\{ i_n \}_{n\in J_k^o}$
satisfying
\begin{equation}
\label{dp1}
\bio{z_n}{z_n^*}  \in
  \BIO{S\left(\wt{\mathcal  P}_k \right)}
  {\left(\lp 1+\e\rp M \right) B\left(\wt{\mathcal  Q}_k \right)}
\end{equation}
and
\begin{equation}
\label{dp2}
\lnm w_{i_n} - z_n \rnm<\delta_k \qquad\text{and}
\qquad\lnm  w^*_{i_n} - z^*_n  \rnm  ~<~ \delta_k  ~+~
\frac{2\delta_k \lp M+\e\rp}{1-2\delta_k}
\end{equation}
for each $n\in J_k^o$.
Towards this,
fix $j\in J_k^o$   and assume that
a biorthogonal system
\begin{equation*}
\bio{z_n}{z_n^*}_{2 + n_k \leq n < j}
\end{equation*}
along with  $\{ i_n \}_{  2 + n_k \leq n < j  } $
have been constructed so that conditions~\eqref{dp1}
and~\eqref{dp2}
hold for $2 + n_k \leq n < j$.
Let
$$
\X_{j} ~\deq~
\wt{\mathcal  P}_k \cap \left[ z^*_n \right]_{2+n_k \leq n < j  }^\pperp
\text{\hskip  25 pt and \hskip  25 pt}
\Y_{j} ~\deq~
\wt{\mathcal  Q}_k \cap \left[ z_n \right]_{2+n_k \leq n < j  }^\perp \ .
$$
Then by Fact~\ref{f:bf} there
exists a natural number $i_j > i_{j-1}$
along with
$z_j\in\X_j$ and $\widetilde{z}_j^*\in\Y_j$
such that
\begin{equation*}
d(w_{i_j},\X_j)
~\leq~
\lnm w_{i_j} - z_j \rnm
~<~
\frac{\delta_k}{M+\e}
\hsp{15}\text{and}\hsp{15}
d(w_{i_j}^*,\Y_j)
~~\leq
\lnm w_{i_j}^* - \widetilde{z}_j^* \rnm
~<~
\delta_k
\end{equation*}
and
\begin{equation*}
\lnm  z_j \rnm ~=~ 1
\hsp{20}\text{and}\hsp{20}
1~\leq~ \lnm \widetilde{z}_j \rnm ~\leq ~M+\e ~.
\end{equation*}
Note that $\widetilde{z}_j^*(z_j)$ need not be equal to $1$ but it
is close to $1$ since
\begin{align}
\begin{split}
\label{e:close}
\lav \widetilde{z}_j^*(z_j)-w_{i_j}^*(w_{i_j}) \rav
~& =~
\lav \widetilde{z}_j^*(z_j)
-(w_{i_j}^*-\widetilde{z}_j^*)(w_{i_j})
-\widetilde{z}_j^*(w_{i_j})\rav
\\
& = ~
\lav \widetilde{z}_j^*(z_j-w_{i_j})-(w_{i_j}^*
-\widetilde{z}_j^*)(w_{i_j}) \rav
\\
& \leq~
\lnm \widetilde{z}_j^*\rnm
\lnm  z_j-w_{i_j} \rnm
+\lnm w_{i_j}^*-\widetilde{z}_j^*\rnm \lnm w_{i_j}\rnm
\\
& <~
\lp M + \e \rp \frac{\delta_k}{M+\e}+\delta_k
~=~
2\delta_k~
\end{split}
\end{align}
and so $1-2\delta_k\leq \widetilde{z}_j^*(z_j)\leq 1+2\delta_k$. Let
\begin{equation*}
 {z}_j^*=\frac{\widetilde{z}_j^*}{\widetilde{z}_j^*(z_j)}
\end{equation*}
so that ${z}_j^*(z_j)=1$.
Now $z_j\in S\left(\tilde{\mathcal{P}}_k\right)$
and $1\leq\| {z}_j^*\|\leq \frac{M+\e}{1-2\delta_k}$
and so  ${z}_j^*\in
\lp 1+2\e\rp {M} B\left(\tilde{\mathcal{Q}}_k\right)$
by ~\eqref{e:Mde}.  Note that by ~\eqref{e:close}
\begin{align*}
\lnm w^*_{i_j} - z_j^* \rnm
~&\leq~
\lnm  w^*_{i_j} - \widetilde{z}_j^* \rnm
~+~
\lnm    \widetilde{z}_j^* - z_j^*  \rnm
\\
~&\leq~
\delta_k ~+~
\lnm  \widetilde{z}_j^* \rnm \
\lav 1 - \frac{1}{\widetilde{z}_j^*(z_j)} \rav
~\leq~
\delta_k ~+~
\lp 1+\e  \rp M \frac{2\delta_k}{1-2\delta_k}~.
\end{align*}
This completes the inductive construction of
$\bio{z_n}{z_n^*}_{n\in J_k^o}$
and $\{ i_n \}_{n\in J_k^o}$.

Now apply the Haar matrix
to $\bio{z_n}{z_n^*}_{n\in J_k}$ to
produce  $\bio{x_n}{x_n^*}_{n\in J_k}$.
With help from the observations in
Remark~\ref{r:hm}, note that
$\bio{x_n}{x_n^*}_{n\in J_k}$
is biorthogonal and is in $\BIO{\mathcal  P_k}{\mathcal  Q_k}$.
Furthermore,
for each $n$ in $J_k$,
\begin{align*}
\lnm x_n \rnm
~&\leq~
2^{{-m_k}/{2}} \   \lnm z_{1+n_k} \rnm ~+~
\left( 1+ \sqrt{2} \right) \max_{j\in J^o_k} \lnm z_{j} \rnm \\
~&\leq~
\e ~+~ \left( 1+ \sqrt{2} \right)
\end{align*}
and
\begin{align*}
\lnm x_n^* \rnm
~&\leq~
2^{{-m_k}/{2}} \   \lnm z_{1+n_k}^* \rnm ~+~
\left( 1+ \sqrt{2} \right) \max_{j\in J^o_k} \lnm z_{j}^* \rnm \\
~&\leq~
\e ~+~ \lp 1+\e \rp M \left( 1+ \sqrt{2} \right).
\end{align*}
If $x^*\in S\left(\X^*\right)$
\begin{align*}
\lav x^{*}\left(x_n\right) \rav
~&\leq~
2^{{-m_k}/{2}} \   \lnm z_{1+n_k} \rnm ~+~
\left( 1+ \sqrt{2} \right)
\max_{j\in J^o_k} \lav x^{*}\left(z_{j}\right) \rav \\
~&\leq~
\delta_k ~+~  \left( 1+ \sqrt{2} \right)
\max_{j\in J^o_k} \left(~ \lav x^{*}\left(w_{i_j}\right) \rav
   +   \delta_k  ~ \right)
\end{align*}
and for each $x^{**}\in S\left(\X^{**}\right)$
\begin{align*}
\lav x^{**} \left( x_n^*\right) \rav~&\leq~
2^{{-m_k}/{2}} \   \lnm z_{1+n_k}^* \rnm ~+~
\left( 1+ \sqrt{2} \right)
\max_{j\in J^o_k} \lav x^{**}\left(z_{j}^*\right) \rav \\
~&\leq~
\delta_k ~+~  \left( 1+ \sqrt{2} \right)
\max_{j\in J^o_k} \left(~ \lav x^{**}\left(w_{i_j}^*\right) \rav
   ~+~   \delta_k
   ~+~ \frac{2\delta_k \lp 1+2\e\rp M}{1-2\delta_k}~ \right) \
\end{align*}
and this simplifies to give us (5).
Thus
$$
  \bio{x_n}{x_n^*}_{n\in  J^p_k \, \cup J_k }
$$
with $\{ i_n \}_{n\in J^{po}_k \cup J^o_k}$
satisfy conditions~(1) through (5).
If $k$ is odd, then
$$
\left[a_h\right]_{h\leq h_k}  \subset \lsp x_n , z_{1+n_k} \rsp_{n\in J^p_k}
        \subset \lsp x_n   \rsp_{n\in J_k^p\cup J_k} \ ,
$$
while if $k$ is even, then
$$
\left[b^*_h\right]_{h\leq h_k} \subset \lsp x_n^* , z_{1+n_k}^* \rsp_{n\in J^p_k}
        \subset \lsp x_n^*   \rsp_{n\in J_k^p\cup J_k} \ .
$$

Clearly the constructed system
$\bio{x_n}{x_n^*}_{n=1}^\infty$,
with  the  blocking $\{ J_k \}_{k=1}^\infty$ of $\mathbb N$ and
the increasing sequence $\{ i_n\}_{n\in \mathbb N^o}$ from $\mathbb N$ satisfy
conditions~(1) through (7).
\end{proof}


\section{Spaces Containing $c_0$}
\label{s:c0}

To motivate the biorthogonal system characterization of spaces containing $c_0$
we recall some well-known facts about such spaces.  We will see that $c_0$
subspaces of $\X$ correspond essentially to weakly unconditionally Cauchy series
in $\X$ so we briefly recall some essential facts about such series.

\begin{definition}
A series $\sum_n x_n$ is called {\it weakly unconditionally
Cauchy} (wuC) if given any permutation $\pi$ of $\N$, the sequence
$\left\{\sum_{k=1}^n x_{\pi (k)}\right\}_n$ is weakly Cauchy.
Equivalently, $\sum_n x_n$ is wuC if and only if for each
$x^*\in\X^*$ we have $\sum_n|x^*(x_n)|<\infty$.
\end{definition}

Bessaga and Pelczynski tied together wuC series and $c_0$ \cite{BP}.

\begin{theorem}
\label{t:c0uvb}\cite{BP}  Let $\X$ be a Banach space.
\begin{enumerate}
\item A basic sequence $\{x_n\}_n$ in $\X$ with $\sum_n x_n$ wuC and $\inf_n \|x_n\|>0$
is equivalent to the unit vector basis of $c_0$.
\item  In order that each wuC series $\sum_n x_n$ in $\X$ be unconditionally convergent
it is both necessary and sufficient that $\X$ contains no copy of $c_0$.
\end{enumerate}
\end{theorem}

Recall the following well-known facts which we will use in this section.

\begin{fact}
\label{f:bas} If $\{x_n\}_n$ is weakly null and $\linf_n\|x_n\|>0$ and $\e>0$,
then $\{x_n\}_n$ has a subsequence which is a basic sequence with basis constant at
most $1+\e$.
\end{fact}

\begin{remark} \label{r:wuCc0}
\begin{enumerate}
\item[(i)] Let $\bio{x_n}{x_n^*}$ be a biorthogonal system with $\sum_n x_n$ wuC and
 $\linf_n \|x_n\|>0$.  If $\{x_{n_k}\}_k$ is any subsequence of $\{x_n\}_n$, then
 $\sum_k x_{n_k}$ is wuC and $\linf_k\|x_{n_k}\|>0$ so Fact \ref{f:bas} tells us
 $\{x_{n_k}\}_k$ has a subsequence $\{x_{n_{k_j}}\}_j$ which is basic and
 $\inf_j\|x_{n_{k_j}}\|>0$. Then by Theorem \ref{t:c0uvb} $\{x_{n_{k_j}}\}_j$ is
 equivalent to the unit vector basis of $c_0$. Thus each subsequence of $\{x_n\}_n$
 has a further subsequence which is equivalent to the unit vector basis of $c_0$.
\item[(ii)] (cf. \cite{D1}) Let $T$ be a bounded linear operator from $c_0$ to $\X$
and \newline$x_n=Te_n$ where $\{e_n\}_n$ is the unit vector basis of $c_0$.  Then
for $x^*\in\X^*$
$$\sum_n|x^*(x_n)|=\sum_n \left|x^*\left(Te_n\right)\right|=\sum_n\left|T^*x^*(e_n)\right|<\infty$$
since $T^*x^*\in\ell_1$.  Thus $\sum_n x_n$ is wuC.  Conversely if
$\sum_n x_n$ is wuC in $\X$, then define $T:c_0\to\X$ by
$T(\{t_n\}_n)=\sum_n t_n x_n$.  Then $T$ is well-defined and has a
closed graph so $T$ is bounded.  So the bounded linear operators
from $c_0$ to $\X$ correspond precisely to the wuC series in $\X$.
\item[(iii)] Let $T:c_0\embed \X$ be an isomorphic embedding and
$\{e_n\}_n$ be the unit vector basis of $c_0$.  Since $T$ is an
embedding there exist constants $C_1$ and $C_2$ such that for any
$(\alpha_n)_n\in c_0$ we have
\[
C_1\|(\alpha_n)_n\|_{c_0}\leq \|T((\alpha_n)_n)\|_\X\leq C_2\|(\alpha_n)_n\|_{c_0}.
\]
Then for each $n\in\N$
\[
C_1=C_1\|e_n\|_{c_0}\leq \|Te_n\|_\X\leq C_2\|e_n\|_{c_0}=C_2
\]
and so $\{Te_n\}_n$ is semi-normalized. By (ii) above, the series $\sum_n Te_n$ is wuC.
\end{enumerate}
\end{remark}

These ideas help us define our $c_0$-biorthogonal system in a very natural way.

\begin{definition} A biorthogonal system $\bio{x_n}{x_n^*}$ in $\BIO{\X}{\X^*}$ is
called a \textit{\cobs} if $\{x_n\}_n$ is normalized and has a subsequence $\{x_{n_j}\}_j$
for which $\sum_j x_{n_j}$ is wuC.
\end{definition}

\begin{theorem}
\label{t:c0it}
The following statements  are equivalent.
\begin{enumerate}
\item $\X$ contains an isomporphic copy of $c_0$.
\item There is a  bounded  \cobs
             in $\BIO{\X}{\X^*}$.
\end{enumerate}
And in the case that $\X$ is separable:
\begin{enumerate}
\item[(3)] There is a  bounded \ft\ \cobs  \newline $\bio{x_n}{x_n^*}\subset\BIO{\X}{\X^*}$.
\end{enumerate}
Furthermore, for each $\e>0$:
if {\rm (2)} holds then the system can be taken to be $(2+\e)$-bounded;
if {\rm (3)} holds then the system
can be taken to be $[2(1+\sqrt2)^2 + \e]$-bounded
and so that  $\lsp x_n^* \rsp $ norms $\X$.
\end{theorem}
\noindent That (2) implies (1) as well as (3) implies (1) follow from Remark \ref{r:wuCc0}.
That (1) implies (2) is Theorem~\ref{t:c0}.  That (1) implies (3) in the separable case
follows from Theorem~\ref{t:ftc0}.

\begin{theorem}
\label{t:c0} If $\X$ contains an isomorphic copy of $c_0$ and $\e>0$, then there exists
a $(2+\e)-$bounded \cobs $\bio{x_n}{x_n^*}\subset\BIO{S(\X)}{\X^*}$.
\end{theorem}

\begin{proof}  Let $T:c_0\embed\X$ be an isomorphic embedding and $\e>0$.  Let $\{e_j\}_j$
be the unit vector basis of $c_0$.  Then by Remark \ref{r:wuCc0} we have $\sum_j Te_j$ is
wuC and $\{Te_j\}_j$ is semi-normalized.  Fact \ref{f:bas} gives us a subsequence
$\{Te_{j_n}\}_n$ of $\{Te_j\}_j$ that is basic with basis constant at most $1+\frac{\e}{2}$.
Let
$$x_n=\frac{Te_{j_n}}{\|Te_{j_n}\|}.$$
Note that $\{x_n\}_n$ is a normalized basic sequence with basis constant at most
$1+\frac{\e}{2}$ and $\sum_n x_n$ is wuC.  We may pick our biorthogonal functionals
accordingly.
\end{proof}
\noindent Notice that the proof of Theorem \ref{t:c0} gives us a
bit more than a \cobs. It gives us a biorthogonal system
$\bio{x_n}{x_n^*}$ with the entire series $\sum_n x_n$ wuC.

To construct a fundamental total biorthogonal system in the separable case we need the
following lemma.

\begin{lemma}\label{l:fincod}  If $Y_0$ is a finite codimensional subspace of $\X^*$ and
$\e>0$, then there is a finite codimensional subspace $X_0$ of $\X$ that is $(2+\e)$-normed
by $Y_0$.
\end{lemma}

\begin{proof}  Let $X_0$ be the pre-annihilator of any finite dimensional subspace of
$\X^*$ that $(1+\e)$-norms the annihilator of $Y_0$.  Then for $f\in S(X_0)$ we have
\begin{align*}
\sup_{y^*\in S(Y_0)}|y^*(f)| &=\inf_{y^{**}\in Y_0^\perp}\|f-y^{**}\| \\
&\geq\inf_{y^{**}\in Y_0^\perp}\text{max}\left[\|f\|-\|y^{**}\|, \sup_{x^*\in S(X_0^\perp)}\left|\lp f-y^{**}\rp(x^*)\right|\right]\\
&\geq \inf_{y^{**}\in Y_0^\perp} \text{max} \left[ 1-\|y^{**}\|, \frac{1}{1+\e}\|y^{**}\|\right]\\
&=\inf_{0\leq t<\infty}\text{max}\left[1-t,\frac{t}{1+\e}\right]\\
&=\frac{1}{2+\e}.
\end{align*}
So $\displaystyle\|f\|\leq (2+\e)\sup_{y^*\in S(Y_0)}|y^*(f)|$ for each $f\in S(X_0)$.
Thus $X_0$ is $(2+\e)$-normed by~$Y_0$.
\end{proof}

The following theorem will give us a fundamental total \cobs in the separable case.

\begin{theorem}
\label{t:ftc0} Suppose $\X$ has a subspace isomorphic to $c_0$.  Let $\e>0$ and
$\bio{a_n}{b_n^*}\subset\BIO{\X}{\X^*}$.  Then there exists a $[2(1+\sqrt{2})^2+\e]$-bounded
\cobs $\bio{x_n}{x_n^*}\subset\BIO{\X}{\X^*}$ with $[a_n]_n\subseteq [x_n]_n$ and $[b_n^*]_n\subseteq [x_n^*]_n$
\end{theorem}

\begin{proof} Without loss of generality,
$[ a_n ]_{n\in\mathbb N}$ and  $[ b^*_n ]_{n\in\mathbb N}$
are each infinite dimensional.  Since $c_0\embed\X$,
by Theorem~\ref{t:c0},
there is a $(2+\e)$-bounded biorthogonal system $\bio{w_n}{w_n^*}$
in $\BIO{S(\X)}{\X^*}$
with $\sum_n w_n$ wuC.  Fix a sequence $\{ \delta_k \}_{k=1}^\infty$
of positive numbers decreasing to zero with $\sum_k \delta_k~<~\infty$. Again we follow the
notation in Definition \ref{d:blk}.  It suffices to find a system
$\bio{x_n}{x_n^*}_{n=1}^\infty$ in $\BIO{\X}{\X^*}$
along with
a blocking $\{ J_k \}_{k=1}^\infty$ of~$\mathbb N$
and
an increasing sequence $\{ i_n\}_{n\in \mathbb N^o}$ from $\mathbb N$,
satisfying
\begin{enumerate}
\item[\rm{(a)}]
 $x_m^*(x_n) = \delta_{mn}$
\item[\rm{(b)}]
$\lnm x_n \rnm ~\leq~ \left(1+\sqrt{2} \right) + \e$
\item[\rm{(c)}]
$\lnm x_n^* \rnm ~\leq~  \lp 2+\e\rp\left(1+\sqrt{2} \right) + \e $
\item[\rm{(d)}]
for each $x^*\in S\left(\X^*\right)$, if $n\in J_k$, then
\newline
$\lav x^{*}\left(x_n\right) \rav ~\leq~ \lp 2+\sqrt{2}\rp\delta_k+\lp 1+\sqrt{2}\rp\text{max}_{j\in J_k^o}\lav x^*(w_{i_j})\rav$
\item[\rm{(e)}]
 $\lsp a_n \rsp_{n=1}^{\infty}
       \subset \lsp x_n \rsp_{n=1}^{\infty}$
\item[\rm{(f)}]
  $\lsp b^*_n \rsp_{n=1}^{\infty}
       \subset \lsp x_n^* \rsp_{n=1}^{\infty} $   \ .
\end{enumerate}

The construction will inductively produce
blocks $\bio{x_n}{x_n^*}_{n\in J_k}$.
Let $x_0$ and $x^*_0$ be the zero vectors.
Fix $k \geq 1$.
Assume that $\{ J_j \}_{0\leq j<k}$
along with
$\bio{x_n}{x_n^*}_{n \in J_k^p}$
and  $\{ i_n \}_{n \in J_k^{po}}$ have been constructed to satisfy conditions~(a) through (d).
Now to construct $J_k$
along with $\bio{x_n}{x_n^*}_{n \in J_k }$
and   $\{ i_n \}_{n  \in J_k^o }$.

Let
$$
 \mathcal P_k ~\deq~
    \left[ x^*_n \right]_{n\in J_k^p }^\pperp
    \text{\hskip 20 pt and \hskip 20 pt}
 \mathcal  Q_k ~\deq~
    \left[ x_n \right]_{n \in J_k^p}^\perp
$$
and
$$
    n_k = \max  J^p_k \ .
$$
The idea is to    find a biorthogonal system
$\bio{z_n}{z_n^*}_{n\in J_k}$ in $\BIO{\mathcal  P_k}{\mathcal
Q_k}$ by first finding  just one pair $\{z_{1+n_k}, z^*_{1+n_k}
\}$ which helps guarantee condition~(e) if $k$ is odd  and
condition~(f) if $k$ is even; however,   $\{z_{1+n_k}, z^*_{1+n_k}
\}$ would not necessarily satisfy conditions~(b) through~(d) so
$J^o_k$ and
$$
    \bio{z_n}{z_n^*}_{n \in J_k^o} \ ,
$$
and $\{ i_n \}_{n\in J^o_k}$ are constructed and then the appropriate Haar matrix is
applied to
$\bio{z_n}{z_n^*}_{n\in J_k}$ to produce $\bio{x_n}{x_n^*}_{n\in J_k}$
so that
$$
  \bio{x_n}{x_n^*}_{n\in  J^p_k \, \cup J_k }
$$
with $\{ i_n \}_{n\in J^{po}_k \cup J^o_k}$
satisfy conditions~(a) through (d).

$\bio{z_{1+n_k}}{z^*_{1+n_k}}$ is constructed
by a standard Gram-Schmidt biorthogonal procedure.
If $k$ is odd, start  in~$\X$. Let
$$
h_k ~=~ \min \left\{ h \colon
     a_h \not\in \lsp x_n \rsp_{n \leq n_k} \right\} \ .
$$
Set
\begin{align*}
z_{1+n_k}
~&=~
a_{h_k} ~-~
   \sum_{n \leq n_k} x_n^*(a_{h_k}) x_n \ ,  \\
\intertext{and for any $y^*_{1+n_k}$ in $\X^*$ such
that   $y^*_{1+n_k}(z_{1+n_k}) \neq 0$,}
z^*_{1+n_k}
~&=~
\frac{
y^*_{1+n_k}  ~-~
   \sum_{n \leq n_k}  y^*_{1+n_k}  (x_n) x_n^*}
{  y^*_{1+n_k}  (z_{1+n_k})} \ .
\end{align*}
If $k$ is even, start  in~$\X^*$. Let
$$
h_{k} ~=~ \min \left\{ h \colon
     b^*_h \not\in \lsp x_n^* \rsp_{n\leq n_k} \right\} \ .
$$
Set
\begin{align*}
z_{1+n_k}^*
~&=~
b^*_{h_k} ~-~
   \sum_{n \leq n_k} b^*_{h_k} (x_n) x_n^* \ ,  \\
\intertext{and, for any $y_{1+n_k}$ in $\X$  such
that   $ z_{1+n_k}^* (y_{1+n_k})  \neq 0$,}
z_{1+n_k}
~&=~
\frac{
y_{1+n_k}  ~-~
   \sum_{n\leq n_k}  x^*_n ( y_{1+n_k} )  x_n}
{ z_{1+n_k}^* (y_{1+n_k}) } \ .
\end{align*}
Clearly  $z_{1+n_k}^* \left( z_{1+n_k} \right) = 1$
and
$$
  z_{1+n_k} ~ \in ~ \mathcal  P_k
  \text{\hskip 25 pt and \hskip 25 pt}
    z_{1+n_k}^* ~ \in ~ \mathcal  Q_k   \ .
$$

Find a natural number  $m_k$ larger than one  so that
$$
2^{{-m_k}/{2}} \
\max\left(~ \lnm z_{1+n_k} \rnm ~,~ \lnm z_{1+n_k}^* \rnm ~\right)
~<~
\min\left(~ \e ~,~ \delta_k ~\right)
$$
and let
$$
  J_k \deq \{ 1 + n_k , \ldots ,  2^{m_k} + n_k \}
\text{\hskip 10 pt and  so \hskip 10 pt}
  J_k^o \deq \{ 2 + n_k , \ldots , 2^{m_k} + n_k \} \ .
$$

Let
$$
\wt{\mathcal  P}_k ~\deq~  \mathcal  P_k \cap \left[ z^*_{1+n_k}\right]^\pperp
\text{\hskip 25 pt and \hskip 25 pt}
\wt{\mathcal  Q}_k ~\deq~ \mathcal  Q_k \cap \left[z_{1+n_k}\right]^\perp \ .
$$
Now we find a biorthogonal system
$\bio{z_n}{z_n^*}_{n\in J_k^o}$
along with $\{ i_n \}_{n\in J_k^o}$
satisfying
\begin{equation}
\label{c01}
\bio{z_n}{z_n^*}  \in
  \BIO{S\left(\wt{\mathcal  P}_k \right)}
  {\left(2+\e \right)B\left(\wt{\mathcal  Q}_k \right)}
\end{equation}
and
\begin{equation}
\label{c02}
\lnm w_{i_n} - z_n \rnm< \delta_k
\end{equation}
for each $n\in J_k^o$.
Towards this,
fix $j\in J_k^o$   and assume that
a biorthogonal system
\begin{equation*}
\bio{z_n}{z_n^*}_{2 + n_k \leq n < j}
\end{equation*}
along with  $\{ i_n \}_{  2 + n_k \leq n < j  } $
have been constructed so that conditions~\eqref{c01}
and~\eqref{c02}
hold for $2 + n_k \leq n < j$.
Let
$$
\X_{j} ~\deq~
\wt{\mathcal  P}_k \cap \left[ z^*_n \right]_{2+n_k \leq n < j  }^\pperp
\text{\hskip  25 pt and \hskip  25 pt}
\Y_{j} ~\deq~
\wt{\mathcal  Q}_k \cap \left[ z_n \right]_{2+n_k \leq n < j  }^\perp \ .
$$
Apply Lemma \ref{l:fincod} with  $Y_0=\Y_j$ to get a finite codimensional subspace
$X_0$ of $\X$ that is $\left(2+\frac{\e}{2}\right)$-normed by $\Y_j$.  Then by
Fact~\ref{f:bf} there
exists a natural number $i_j > i_{j-1}$
along with
$z_j\in S\left(\X_j\cap X_0\right)$ such that
\begin{equation*}
d(w_{i_j},\X_j\cap X_0)
~\leq~
\lnm z_j-w_{i_j} \rnm
~<~
\delta_k
\end{equation*}
Since $X_0$ is $\left(2+\frac{\e}{2}\right)$-normed by $\Y_j$ there is
$\widetilde{z}_j^*\in S(\Y_j)$ such that
$$\frac{1}{2+\e}\leq\widetilde{z}_j^*(z_j).$$
Let
$$z_j^*=\frac{1}{\widetilde{z}_j^*(z_j)}\widetilde{z}_j^*$$
so that $z_j^*(z_j)=1$ and note that
$$\|z_j^*\|=\frac{1}{\widetilde{z}_j^*(z_j)}\|\widetilde{z}_j^*\|\leq 2+\e.$$
This completes the inductive construction of
$\bio{z_n}{z_n^*}_{n\in J_k^o}$
and $\{ i_n \}_{n\in J_k^o}$.

Now apply the Haar matrix
to $\bio{z_n}{z_n^*}_{n\in J_k}$ to
produce  $\bio{x_n}{x_n^*}_{n\in J_k}$.
With help from the observations in
Remark~\ref{r:hm}, note that
$\bio{x_n}{x_n^*}_{n\in J_k}$
is biorthogonal and is in $\BIO{\mathcal  P_k}{\mathcal  Q_k}$.
Furthermore,
for each $n$ in $J_k$,
\begin{align*}
\lnm x_n \rnm
~&\leq~
2^{{-m_k}/{2}} \   \lnm z_{1+n_k} \rnm ~+~
\left( 1+ \sqrt{2} \right) \max_{j\in J^o_k} \lnm z_{j} \rnm \\
~&\leq~
\e ~+~ \left( 1+ \sqrt{2} \right)
\end{align*}
and
\begin{align*}
\lnm x_n^* \rnm
~&\leq~
2^{{-m_k}/{2}} \   \lnm z_{1+n_k}^* \rnm ~+~
\left( 1+ \sqrt{2} \right) \max_{j\in J^o_k} \lnm z_{j}^* \rnm \\
~&\leq~
\e ~+~ \lp 2+\e \rp  \left( 1+ \sqrt{2} \right).
\end{align*}
If $x^*\in S\left(\X^*\right)$
\begin{align*}
\lav x^{*}\left(x_n\right) \rav
~&\leq~ 2^{-m_k/2}\|z_{1+n_k}\| +\lp 1+\sqrt{2}\rp\text{max}_{j\in J_k^o}|x^*(z_j)|\\
&\leq \delta_k+\lp 1+\sqrt{2}\rp\text{max}_{j\in J_k^o}\lp\left|x^*(z_j-w_{i_j})\right|+\left|x^*(w_{i_j})\right|\rp\\
&\leq\delta_k+\lp 1+\sqrt{2}\rp\text{max}_{j\in J_k^o}\lp\delta_k+\left|x^*(w_{i_j})\right|\rp\\
&=\lp 2+\sqrt{2}\rp\delta_k+\lp 1+\sqrt{2}\rp\text{max}_{j\in J_k^o}\left|x^*(w_{i_j})\right|.
\end{align*}

Thus
$$
  \bio{x_n}{x_n^*}_{n\in  J^p_k \, \cup J_k }
$$
with $\{ i_n \}_{n\in J^{po}_k \cup J^o_k}$
satisfy conditions~(a) through (d).
If $k$ is odd, then
$$
\left[a_h\right]_{h\leq h_k}  \subset \lsp x_n , z_{1+n_k} \rsp_{n\in J^p_k}
        \subset \lsp x_n   \rsp_{n\in J_k^p\cup J_k} \ ,
$$
while if $k$ is even, then
$$
\left[b^*_h\right]_{h\leq h_k} \subset \lsp x_n^* , z_{1+n_k}^* \rsp_{n\in J^p_k}
        \subset \lsp x_n^*   \rsp_{n\in J_k^p\cup J_k} \ .
$$

Clearly the constructed system
$\bio{x_n}{x_n^*}_{n=1}^\infty$,
with  the  blocking $\{ J_k \}_{k=1}^\infty$ of $\mathbb N$ and
the increasing sequence $\{ i_n\}_{n\in \mathbb N^o}$ from $\mathbb N$ satisfy
conditions~(a) through (f).

Note that condition (d) tells us that if for each $k\in\N$ we pick any $n_k\in J_k$, then
for $x^*\in S(\X^*)$ we have
\begin{align*} \sum_k \lav x^*(x_{n_k})\rav &\leq \lp 2+\sqrt{2}\rp\sum_k\delta_k+\lp 1+\sqrt{2}\rp\sum_k\text{max}_{j\in J_k^o}\lav x^*(w_{i_j})\rav\\
&\leq \lp 2+\sqrt{2}\rp\sum_k\delta_k+\lp 1+\sqrt{2}\rp\sum_j \lav x^*(w_{i_j})\rav<\infty
\end{align*}
So $\sum_k x_{n_k}$ is wuC.

\end{proof}


\section{Piecing it all Together}
\label{s:pieces}

Inspired by Theorem \ref{t:Bs} we might try to combine Theorems \ref{t:dpit} and
\ref{t:c0it} with the Dilworth, Girardi, Johnson $\ell_1$ result (Theorem \ref{t:elle})
to get the following theorem giving the existence biorthogonal systems in any Banach space.

\begin{pottheorem}
\label{t:Bsbs1} For any given infinite dimensional Banach space $\X$ there exists a
bounded biorthogonal system $\bio{x_n}{x_n^*}$ that is one of the following three types:
\begin{enumerate}
\item a \cobs
\item a $wc_0^*$-stable biorthogonal system
\item a DP-biorthogonal system.
\end{enumerate}
\end{pottheorem}

\noindent However, this does not follow directly from the previous results.  The
trouble lies in part (3).   Theorem~\ref{t:Bs} guarantees us that if $\X$ contains
no isomorphic copies of $c_0$ or $\ell_1$, then there is a {\it subspace} (say $\Y$)
of $\X$ that fails DP.  So from Theorem \ref{t:dpit} we get a DP-biorthogonal system
$\bio{y_n}{y_n^*}$ in $\BIO{\Y}{\Y^*}$.  Since $\{y_n\}_n$ is weakly null in $\Y$ it
is also weakly null in $\X$.  Unfortunately the fact that $\{y_n^*\}_n$ is weakly null
in $\Y^*$ does not necessarily tell us that if we extend each $y_n^*$ to $x_n^*\in\X^*$,
then $\{x_n^*\}_n$ is weakly null in $\X^*$.  Another way to see that part (3) is not
correct is to notice that DP does not necessarily pass to closed subspaces.  Since it
is a $C(K)$ space, $\ell_\infty$ has DP; however $\ell_2$ does not have DP.  So if part
(3) were correct it would say that $\Y$ failing DP implies $\X$ fails DP, which is false.
We recall the following related property.

\begin{definition} A Banach space $\X$ has the {\it hereditary Dunford-Pettis property}
($\text{DP}_h$) if every closed subspace of $\X$ has the Dunford-Pettis property.
\end{definition}

For detailed discussions of \DPh see \cite{CG,Cem,D2}. In 1987 Cembranos gave the
following useful characterization of \DPh.

\begin{theorem}\label{t:DPh}\cite{Cem} A Banach space $\X$ has \DPh if and only if
every normalized weakly null sequence in $\X$ has a subsequence which is equivalent
to the unit vector basis of $c_0$.
\end{theorem}

\noindent In 1989 Knaust and Odell \cite{KO} gave a quantitative improvement of this
result by showing that the equivalence is uniform for all normalized weakly null sequences.
Using the hereditary Dunford-Pettis property we can restate Theorem \ref{t:Bs}.

\begin{restatement}\label{t:Bs1} Every infinite dimensional Banach space, $\X$,
contains a subspace isomorphic to $c_0$, a subspace isomorphic to $\ell_1$ or $\X$ fails \DPh.
\end{restatement}

\noindent In light of this restatement we see that a biorthogonal system characterization
of \DPh is in order.  Theorem \ref{t:DPh} will give it to us.

\begin{definition} A biorthogonal system $\bio{x_n}{x_n^*}$ in $\BIO{\X}{\X^*}$ is
called a {\it \DPhbs} if $\{x_n\}_n$ is semi-normalized, weakly null and for any
subsequence $\{x_{n_j}\}_j$ the series $\sum_j x_{n_j}$ is not wuC.
\end{definition}

\begin{theorem}\label{t:DPhbs} A Banach space $\X$ fails \DPh if and only if for
each $\e>0$ there is a $(2+\e)$-bounded \DPhbs $\bio{x_n}{x_n^*}$ in $\BIO{S(\X)}{\X^*}$.
\end{theorem}

\begin{proof} ($\Rightarrow$) Suppose $\X$ fails \DPh and $\e>0$.  Then Theorem
\ref{t:DPh} gives us a normalized weakly null sequence $\{x_n\}_n$ with no subsequence
equivalent to the unit vector basis of $c_0$.  Without loss of generality $\{x_n\}_n$
is a basic sequence with basis constant at most $2+\e$.  Now if for some subsequence $\{x_{n_j}\}_j$
we have $\sum_j x_{n_j}$ wuC then Theorem \ref{t:c0uvb} tells us that $\{x_{n_j}\}_j$
is equivalent to the unit vector basis of $c_0$, which is a contradiction.  Since
$\{x_n\}_n$ is basic with basis constant at most $2+\e$, we may pick a sequence of
biorthogonal functionals $\{x_n^*\}_n\subset\lp2+\e\rp B(\X^*)$.

($\Leftarrow$) Suppose there exists such a biorthogonal system $\bio{x_n}{x_n^*}$.
If $\X$ has \DPh then Theorem \ref{t:DPh} gives us a subsequence $\{x_{n_j}\}_j$ of
$\{x_n\}_n$ that is equivalent to the unit vector basis of $c_0$.  But then we would
have $\sum_j x_{n_j}$ wuC which is a contradiction.
\end{proof}

Finally, putting this together with Theorems \ref{t:dpit} and \ref{t:c0it} and the
Dilworth, Girardi, Johnson $\ell_1$ result we get a correct theorem..

\begin{theorem}\label{t:Bsbs2} For any given infinite dimensional Banach space $\X$
there exists a bounded biorthogonal system $\bio{x_n}{x_n^*}$ that is one of the
following three types:
\begin{enumerate}
\item a \cobs
\item a $wc_0^*$-stable biorthogonal system
\item a $\text{DP}_h$-biorthogonal system.
\end{enumerate}
\end{theorem}

Note that this theorem confirms the importance of $c_0$ in
infinite dimensional Banach spaces.  The presence of a \cobs
$\bio{x_n}{x_n^*}$ in $\X$ gives us a part of $\X$ which is
particularly $c_0$-rich in the sense that $[x_n]$ is isomorphic to
$c_0$ by design and, of course, the same is true for any
subsequence $\{x_{n_j}\}_{j=1}^\infty$.  On the other hand, the
existence of a $\text{DP}_h$-biorthogonal system
$\bio{x_n}{x_n^*}$ in $\X$ would signify a part of $\X$ is
completely lacking in $c_0$ subspaces.  In particular, $[x_n]$ is
not isomorphic to $c_0$ and the same is true for any subsequence
$\{x_{n_j}\}_{j=1}^\infty$ since $\sum_n x_{n_j}$ is not wuC. In
the third case if $\X$ has a $wc_0^*$-stable biorthogonal system
$\bio{x_n}{x_n^*}$, then $[x_n]$ is not isomorphic to $c_0$ since
the proof in \cite{DGJ} yields that $[x_n]\approx \ell_1$.

It would be interesting to see what this interpretation of Theorem
\ref{t:Bsbs2} yields in terms of other properties and structures
that have been characterized using $c_0$.  For instance, can we
say anything about the existence of spreading models or nice
(resp. not very nice) operators on the space?

\bibliographystyle{amsplain}
\bibliography{biblio}
\end{document}